\documentclass[11pt, reqno]{amsart}% cls. for J. Nonlinear Sci. Appl

\usepackage{amsmath, amsthm, amscd, amsfonts, amssymb, graphicx, color, enumerate}
\RequirePackage[backref,colorlinks=true]{hyperref}
\newcommand{\R}{{\mathbb{R}}}
 \newcommand{\N}{{\mathbb{N}}}
 \newcommand{\inter}{\mathopen [0,1\mathclose]}
\newcommand{\interu}{\mathopen [u\sp{-}(\lambda),
u\sp{+}(\lambda)\mathclose]}
\newcommand{\interv}{\mathopen [v\sp{-}(\lambda),
v\sp{+}(\lambda)\mathclose]}

 \newcommand{\m}{\medskip}
 
 \newcommand{\cj}[2]{\left \{\, #1\, : \, #2\, \right\}}
 
 \newcommand{\n}{\noindent}

\newcommand{\fn}{\mathbb{E}\sp{1}}

 \DeclareMathOperator{\cl}{cl}

\begin{document}

\title[Jackson-type approximation for fuzzy-valued functions ]{Jackson-type approximation for fuzzy-valued functions by means of trapezoidal functions}

%    Remove any unused author tags.

%    author one information
\author{Juan J. Font}
\address{Institut Universitari
de Matem\`{a}tiques i Aplicacions de Castell\'{o} (IMAC),
Universitat Jaume I, Avda. Sos Baynat s/n, 12071, Castell\'{o}
(Spain).}
\thanks{Research supported by grant PID2019-106529GB-I00 and grant
UJI-B2019-08}
%\curraddr{}
\email{font@uji.es}

%    author two information
\author{Sergio Macario}

\address{Institut Universitari
de Matem\`{a}tiques i Aplicacions de Castell\'{o} (IMAC),
Universitat Jaume I, Avda. Sos Baynat s/n, 12071, Castell\'{o}
(Spain).}
%\thanks{The second author is supported by grant  UJI-B2019-08}
%\curraddr{}
\email{macario@uji.es}



\date{\today}

% {\tt
%\begin{center}Communicated by ??? \end{center}}}
%

%The name will be entered by editor

%\cortext[c1]{Corresponding author}
%\setcounter{page}{1}
%\vol{? (201?)} \pages{1--?}
%% The correct information will be entered by the editor
%
%\recivedat{received date}
%% The correct dates will be entered by the editor

%\authors{J.J. Font, S. Macario}

%\doi{\href{doi.org/10.22436}{10.22436/jnsa.010.01.01}}
% The correct doi will be entered by the editor

\begin{abstract}
In this paper we provide new several Jackson-type approximations results for continuous fuzzy-number-valued functions which improve several previous ones. We use alternative techniques adapted from Interval Analysis which rely on the gH-difference (which might not exist)  and the generalized difference (which might lack the cancellation law ) of fuzzy numbers.
\end{abstract}

\maketitle

%Theorem-like structures
%If you need new environments, define them here with the command \newtheorem{envname}{caption}.

%Some environments such as definitions, theorems, lemmas or examples, have defined  in the list below.

\newtheorem{theorem}{Theorem}[section]
\newtheorem{lemma}[theorem]{Lemma}
\newtheorem{proposition}[theorem]{Proposition}
\newtheorem{corollary}[theorem]{Corollary}
\newtheorem{question}[theorem]{Question}
\newtheorem{example}[theorem]{Example}

\theoremstyle{definition}
\newtheorem{definition}[theorem]{Definition}
\newtheorem{algorithm}[theorem]{Algorithm}
\newtheorem{conclusion}[theorem]{Conclusion}
\newtheorem{problem}[theorem]{Problem}

\theoremstyle{remark}
\newtheorem{remark}[theorem]{Remark}
\numberwithin{equation}{section}

\section{Introduction}

One of the main goals in fuzzy approximation theory is to approximate a
continuous fuzzy-valued function defined on a compact space by some simpler functions such as (trigonometric) polynomials, rational functions,  piecewise linear functions, etc. (see \cite{anas}).

The most common scheme is to find uniform approximation theorems of Weierstrass type 
(\cite{gal:93}, \cite{FSS:17}, \cite{syau}, \cite{PANA}). But, in order to set a quantitative estimation of the degree of approximation, the error bounds are given in terms of the fuzzy modulus of continuity (see e.g. \cite{gal:94}, 
\cite{BedeGal:2004}), the so-called Jackson-type estimates (\cite{jackson}).

We will focus on the following scheme: let $\fn$ denote the set of fuzzy numbers (see Section 2) and let us consider the family of functions
$$\mathfrak{G}_n:=\left\{\sum_{k=0}^{n} \phi_k(x)\cdot u_k\ :\  \phi_k(x) \text{ real-valued continuous functions},  u_k\in\fn\right\}$$
and, for a given fuzzy-valued continuous function $f\in C([a,b],\fn)$,  denote the corresponding approximation error
$$E_{n,f}:=\inf\{D(f,g) \ :\ g\in\mathfrak{G}_n\},$$
where
$D$ is the supremum metric on
$C([a,b],\fn)$.

In this paper we are interested in giving an estimate of this approximation
 in the form
 $$E_{n,f}\leq C\omega^{\mathcal{F}}\left(f,\frac{1}{n}\right)$$
 where $\omega^{\mathcal{F}}\left(f,\frac{1}{n}\right)$ stands for the fuzzy modulus of continuity for such functions (defined below), and $C$ is a constant independent of $n$ and $f$.

Several such schemes have   been already developed. For instance,
when the coefficients $\phi_k$ are  polynomials,
S. Gal proved in \cite{Gal:2000} that, given $f\in C([0,1],\fn)$, then
$$D(f,B_n^{\mathcal{F}}(f))\leq \dfrac{3}{2}\omega^{\mathcal{F}}\left(f,\frac{1}{\sqrt{n}}\right)$$
where
$$B_n^{\mathcal{F}}(f) :=\sum_{k=0}^n \binom{n}{k} x^k(1-x)^{n-k} f\left(\frac{k}{n}\right)$$
are called fuzzy polynomials of Bernstein type. To improve the rate of convergence from $\frac{1}{\sqrt{n}}$ to $\frac{1}{n}$, Bede and Gal (\cite{BedeGal:2004}) provided a Jackson-type estimate using Szabados-type polynomials but with an unexplicit constant~$C$. Also they gave an estimation for approximating $2\pi$-periodic continuous fuzzy-valued functions by sequences of fuzzy trigonometric polynomials, completing a previous work of Anastassiou and Gal (\cite{AG:01}).
Also in \cite[Lemma 3.4]{BedeGal:2004}, following the ideas of Szabados (\cite{sza}), they obtained a Jackson-type estimate for rational functions using Shepard-like operators.
%All these results can be found in Chapter 5 and 6 in \cite{anas} which are devoted to present the results on trigonometrical and polynomial fuzzy approximation.
Other recent papers in the topic  give similar estimates of the approximation using a generalization of Bernstein-type polynomials, the so-called post-quantum Bernstein polynomials   (\cite{ozkan})  or the fuzzy $(p,q)$-Bernstein-Chlodowsky operators (\cite{ozkan2}) or fuzzy Berstein-type rational functions (\cite{ozkan3}).

The technique used in all those papers requires that, in order to compute $D(f,\sum \phi_k(x)\cdot u_k)$,  the sum of the coefficients $\sum_{k=0}^{n} \phi_k(x)=1$.
%Our last theorem in this paper follows this approach, obtaining the following approximation result:
%$$D\left(f,\sum_{k=0}^{n}\psi_k(x)\cdot f\left(\frac{k}{n}\right)\right)\leq 3\omega^{\mathcal{F}}\left(f, \frac{1}{n}\right).$$
%Here the coefficients $\psi_j(x)$ are trapezoidal functions defined from $[0,1]$ into $[0,1]$.

In this paper we will try another approach which aligns with the ideas developed in \cite{FontMac:2022},
where the authors provided a Jackson-type approximation for continuous interval-valued functions defined on an interval $[a,b]$, combining techniques used by Prolla (\cite{prola}) to give a constructive proof of the Weierstrass-Stone Theorem and D. Chen (\cite{Chen}) to estimate the approximation error by neural networks.
%The main result in \cite{FontMac:2022} states that, given a continuous interval-valued function $f$
%and for each $\epsilon>0$ we can find a  function
%$$g(x)=\sum_{j=0}^{n}\psi_j(x)A_j$$
%where $\psi_j(x)$ are continuous functions from $[0,1]$ into itself and $A_j$ are compact intervals of the real line,
%and with the property that
%$$D(f,g)\leq 2\omega^{\mathcal{K}}\left(f, \frac{1}{n}\right)+\epsilon.$$
%Here $\omega^{\mathcal{K}}\left(f, \frac{1}{n}\right)$ stands for the modulus of continuity for interval-valued continuous functions.
We intend to export such technique to  provide  similar results for continuous fuzzy-valued functions in a 
essentially different way from the mentioned papers since it is not requiring that $\sum_{j=0}^{n}\psi_j(x)=1$ (in fact, it will be greater than one).

It is well known that a fuzzy number can be identified as a family of nested compact intervals, the so called  $\alpha$-cuts. So it is expected that our previous research for continuous interval-valued functions can be applied to get a sharp result, that is, given now $f\in C([0,1],\fn)$, our objective is to find a fuzzy-valued function like
\begin{equation}\label{eq:gg}
g(x)=\sum_{j=0}^{n}\psi_j(x)u_j
\end{equation}
where now $u_j\in\fn$ with
$$D(f,g)\leq C\omega^{\mathcal{F}}\left(f, \frac{1}{n}\right)+\epsilon.$$

However, the main drawback in the fuzzy-valued context is that Chen's techniques require some sort of cancelation law, since such techniques use $u_j$ in Equation~\eqref{eq:gg} as a difference between two fuzzy numbers. This fact was not a problem for interval-valued functions because the generalized Hukuhara difference, gH-difference for short,  has the cancellation law (\cite{Stefanini:2008}).  So,  we would expect to obtain a similar approximation by replacing the compact intervals by fuzzy numbers and using the gH-difference of fuzzy numbers. Nevertheless,  such a substraction does not always provide a fuzzy number. So we have to restrict our attention to some kind of fuzzy-valued functions for which our results still work
 (Theorem~\ref{th:gH-fuzzy degree}).

In order to broaden the types of functions to which we can apply our technique, we replace the gH-difference of fuzzy numbers with their generalized difference (\cite{BedeStef:2013}).  In this case, our scope is larger than in  Theorem~\ref{th:gH-fuzzy degree}, but our bound is not sharp as it grows with $n$ (Theorem~\ref{th:main-fuzzy-degree}) due to the fact that the generalized difference of fuzzy number lacks the cancelation law. So it seems that Chen's technique is not the best one for fuzzy-valued functions.

Finally, Prolla's ideas still allow us to construct a family of piecewise linear functions which provides a good approximation bound, although a bit worse than that of Theorem~\ref{th:gH-fuzzy degree}, but which is valid for any continuous fuzzy-valued function (Theorem~\ref{th:degree-approx-fuzzy}).

% Namely,
%$$D\left(f,\sum_{k=0}^{n}\psi_k(x)\cdot f\left(\frac{k}{n}\right)\right)\leq 3\omega^{\mathcal{F}}\left(f, \frac{1}{n}\right)$$
%where now the coefficients $\psi_j(x)$ are  piecewise linear  functions (see Proposition~\ref{th:sum2}) and $\sum_{k=0}^{n}\psi_k(x)=1$.

\section{Preliminaries}

Let $F(\R)$ denote the family of all fuzzy subsets on the real
numbers $\R$, that is, mappings $u:\R\rightarrow [0,1]$.
For $u\in F(\R)$ and $\lambda\in\inter$, the $\lambda$-level set
of $u$ is defined by
$$
[u]\sb{\lambda}:=\cj{x\in\R}{u(x)\geq \lambda}, \quad
\lambda\in ]0,1] ,\quad [u]\sb{0}:=\cl\sb{\R}\cj{x\in\R}{u(x)>0}.
$$

The fuzzy number space $\fn$ is the set of elements $u$ of
$F(\R)$ satisfying the following properties:

\begin{enumerate}[(1)]
 \item $u$ is normal, i.e., there exists an $x\sb{0}\in\R$ with $u(x\sb{0})=1$;

\item $u$ is convex, i.e., $u(\lambda x + (1-\lambda)y)\geq \min
\left\{u(x),u(y)\right\}$  for all $x,y\in \R, \lambda\in\inter$;

\item $u$ is upper-semicontinuous;

\item $[u]\sb{0}$ is a compact set in $\R$.
\end{enumerate}

Notice that if $u\in\fn$, then the $\lambda$-level set
$[u]\sb{\lambda}$ of $u$ is a member of $\mathcal{K_C}$ for each
$\lambda\in\inter$. We denote $[u]\sb{\lambda}=\interu$.  Every real
number $r$ can be
identified with the fuzzy number $\tilde{r}$ defined as
\[
\tilde{r}(t):=\left\{
\begin{array}{lc}
1 & \text{if } t=r,\m\\
0 & \text{if } t\not= r.
\end{array}\right.
\]
We can now state the characterization of fuzzy numbers provided by Goetschel and Voxman (\cite{GV}):

\begin{theorem}\label{GW1}
Let $u\in\fn$ and $[u]\sb{\lambda}=\interu$, $\lambda\in\inter$.
Then the pair of functions $u\sp{-}(\lambda)$ and
$u\sp{+}(\lambda)$ has the following properties:

\begin{enumerate}[(i)]
 \item $u\sp{-}(\lambda)$ is a bounded left continuous nondecreasing function
on $\mathopen (0,1\mathclose]$;

\item  $u\sp{+}(\lambda)$ is a bounded left continuous nonincreasing function
on $\mathopen (0,1\mathclose]$;

\item $u\sp{-}(\lambda)$ and $u\sp{+}(\lambda)$ are right continuous at
$\lambda=0$;

\item $u\sp{-}(1)\leq u\sp{+}(1)$.
\end{enumerate}

\n Conversely, if a pair of functions $\alpha (\lambda)$ and
$\beta(\lambda)$ satisfy the above conditions (i)-(iv), then there
exists a unique $u\in\fn$ such that
$[u]\sb{\lambda}=\mathopen[\alpha(\lambda),\beta(\lambda)\mathclose]$
for each $\lambda\in\inter$.
\end{theorem}

Given $u,v\in \fn$ and $k\in \mathbb{R}$, we can define $u+v$ and $ku$ level-wise as $[u+v]_{\lambda}:=\interu+\interv$ and $[ku]_{\lambda}:=k\interu$.
It is well-known that $\fn$ endowed with these two natural operations is not a vector space.
In fact, there is no inverse for the sum operation.
On the other hand, we can endow $\fn$ with the following metric:

\begin{definition}{\rm \cite{GV,DK2}}\label{GW2}
 For $u, v\in\fn$, we can define

$$
d\sb{\infty}(u,v):=\sup\sb{\lambda\in\inter}\max\left\{|u\sp{-}(\lambda)-v\sp{
- } (\lambda)|, |u\sp{+}(\lambda)-v\sp{+} (\lambda)|\right \}.
$$
It is worth to remark that

$$
d\sb{\infty}(u,v):=\sup\sb{\lambda\in\inter}\left\{d_H([u]_{\lambda},[v]_{\lambda}) \right\}.
$$
where $d_H$ stands for the Hausdorff metric on the space of all closed intervals.
\end{definition}

$d\sb{\infty}$ is called the supremum
metric on $\fn$, and $(\fn , d\sb{\infty})$ is well-known to be a complete metric
space.
Notice that, by the definition of $d\sb{\infty}$, $\R$
endowed with the euclidean topology can be topologically identified
with the closed subspace $\tilde{R}=\cj{\tilde{x}}{x\in\R}$ of
$(\fn , d\sb{\infty})$ where
$\tilde{x}\sp{+}(\lambda)=\tilde{x}\sp{-}(\lambda)=x$ for all
$\lambda\in\mathopen[0,1\mathclose]$. As a metric space, we shall always
consider $\fn$ equipped with the metric $d\sb{\infty}$.

The following Lemmas are the key to find the coefficients in our approximation function.
\begin{lemma}\cite[Lemma 2]{Je}\label{Jewett} Let $0\leq a < b \leq 1$ and $0 < \epsilon < \frac{1}{2}$. There exists a polynomial $p(x)= (1 - x^m)^n$ such that
\begin{enumerate}
\item $p(x)> 1 - \epsilon$ for all $0 \le x \le a$,
\item $p(x)< \epsilon$ for all $b \le x \le 1$.
\end{enumerate}
\end{lemma}

The lemma above will be used throughout this paper in the following form:
\begin{lemma}\label{lem:Je2}
Let $n$ be positive integer number and set $a_j=\frac{j}{n}$ for $j=0,1,\ldots, n$. Given $0<\epsilon<\frac{1}{2}$ and $0<\delta<\frac{1}{2n}$, we can find $n+1$ continuous functions $\psi_j:[0,1]\rightarrow[0,1]$ such that
\begin{enumerate}
\item $\psi_j(x)<\epsilon$, for $x\geq a_j+\delta$, $j=0,\ldots, n-1$;

 \item $1-\psi_j(x)<\epsilon $, for $x\leq a_j-\delta$, $j=1,\ldots, n$.

% \item $\sum_{j=0}^{n}\psi_j(x)\leq 1$, for all $x\in[0,1]$.
\end{enumerate}

\end{lemma}
\begin{proof}
Apply Lemma~\ref{Jewett} to $a_j-\delta<a_j+\delta$ $j=1,\ldots, n-1$ (or $a_0<a_0+\delta$ for $j=0$ and $a_n-\delta<a_n$ for $j=n$). Therefore,  we are able to find $n+1$ continuous functions $\psi_j: [0,1]\rightarrow[0,1]$ satisfying the requested conditions.
\end{proof}

Let us finally recall some related results concerning interval-valued functions.

Let $\mathcal{K_{C}}$ denote the set of all compact intervals of the real line $\mathbb{R}$; that is
$$\mathcal{K_{C}}=\{[a,b]\ : a,b\in\mathbb{R},\ a\leq b\}.$$
where $[a,a]$ stands for the singleton $\{a\}$. We denote by $len(A)$ the diameter or length of the interval $A$; that is, $len(A)=a^+-a^-$, when $A=[a^-,a^+]$.

Let  us denote the family
$$\mathfrak{T}_n=\left\{\sum_{i=0}^{n} \psi_i(x)A_i\ :\ \psi_i \in C([0,1],[0,1]); \ A_i\in\mathcal{K_C},\ i=0,1,2,\ldots n\right\}.$$
The authors proved in \cite{FontMac:2022} that  the set $\mathfrak{T}=\cup_{n\in\N} \mathfrak{T}_n$ is dense in $C([0,1],\mathcal{K_C})$.
Define the  approximation error between
$\mathfrak{T}_n$ and a continuous function $f\in C([0,1],\mathcal{K_C})$   by
$$
E_{n,f}:=\inf_{g\in \mathfrak{T}_n}D(f,g)
$$
where $D(f,g):=\sup\{d_{H}(f(x),g(x))\ :\ x\in [0,1]\}$ is the supremum metric on $C([0,1],\mathcal{K_C})$.
Let us consider the  modulus of continuity of  $f$ defined as
\[
\omega^{\mathcal{K}}(f,\delta):=\sup\{ d_H(f(x),f(y))\; :\; x,y\in[0,1],\  |x-y|<\delta\}.
\]

Relying on an idea which we can trace back to Chen (\cite{Chen}) and further developed by Hong and Hahm (\cite{HH:2002}, \cite{HH:2016}),  the authors provided in \cite{FontMac:2022}, using the generalized Hukuhara difference for intervals,  a Jackson type approximation for interval-valued continuous functions defined on $[0,1]$  as follows.

\begin{theorem}\cite{FontMac:2022}\label{th:IV-degree1}
Let $f\in C([0,1],\mathcal{K_C})$ satisfying

\begin{equation}\label{eq:length1}
len(f(x))\geq len(f(y)),\quad \text{when $x\leq y$.}
\end{equation}

or

\begin{equation}\label{eq:length2}
len(f(x))\leq len(f(y)),\quad \text{when $x\leq y$.}
\end{equation}

\medskip
Then, for every $n\in\mathbb{N}$,

\[E_{n,f} \leq 2\omega^{\mathcal{K}}\left(f,\dfrac{1}{n}\right).\]

\end{theorem}

In fact, from the proof of that theorem, we can, using Lemma~\ref{lem:Je2}, get the following propositions.

\begin{proposition}{\cite[Theorem 4.3 ]{FontMac:2022}}\label{prop:IV-degree2}
Let $f\in C([0,1],\mathcal{K_C})$ satisfying equation~\eqref{eq:length1} and let $M>1$ be such that
$d_H(f(x),f(y))\leq 2M$, for all $x,y\in[0,1]$.
Let $\varepsilon>0$, $n\in \mathbb{N}$ and $0<\delta<\frac{1}{2n}$.
Let $a_j=\frac{j}{n}$ be given and take $\epsilon'<\frac{\varepsilon}{2(n+1)M}$. Let   $\psi_j:[0,1]\rightarrow [0,1]$ be the  $n+1$ continuous functions given by Lemma~\ref{lem:Je2} for such $\epsilon'$ and $\delta$.

Then,  the function
$$g(x):=A_{n}+\sum_{j=0}^{n-1} \psi_j (x)\Big(A_{j}\circleddash A_{j+1}\Big)$$
where $A_j=f(a_j)$, $j=0,1,\ldots, n$, verifies that, for all $x\in[0,1]$,
$$
d_H(g(x),f(x))\leq 2\omega^{\mathcal{K}}\left(f,\frac{1}{n}\right) +\varepsilon.
$$

\end{proposition}
\begin{proposition}{\cite[Theorem 4.4 ]{FontMac:2022}}\label{prop:IV-degree3}
Let $f\in C([0,1],\mathcal{K_C})$ satisfying equation~\eqref{eq:length2} and let $M>1$ be such that
$d_H(f(x),f(y))\leq 2M$, for all $x,y\in[0,1]$.
Let $\varepsilon>0$, $n\in \mathbb{N}$ and $0<\delta<\frac{1}{2n}$.
Let $a_j=\frac{j}{n}$ be given and take $\epsilon'<\frac{\varepsilon}{2(n+1)M}$. Let   $\psi_j:[0,1]\rightarrow [0,1]$ be the  $n+1$ continuous functions given by Lemma~\ref{lem:Je2} for such $\epsilon'$ and $\delta$.

Then,  the function
$$g(x):=A_{0}+\sum_{j=0}^{n-1} \phi_j (x)\Big(A_{j+1}\circleddash A_{j}\Big)$$
where $A_j=f(a_j)$ and $\phi_j=1-\psi_j$, $j=0,1,\ldots, n$, verifies that, for all $x\in[0,1]$,
$$
d_H(g(x),f(x))\leq 2\omega^{\mathcal{K}}\left(f,\frac{1}{n}\right) +\varepsilon.
$$
\end{proposition}

\section{Approximation of continuous fuzzy-valued functions}

We will present similar results as those mentioned in the section above but  for functions $f\in C([0,1],\fn)$, endowed with the distance
$$D(f,g):=\sup_{x\in[0,1]} d_{\infty}(f(x),g(x)).$$

The generalized Hukuhara difference for fuzzy numbers was defined by  Stefanini in
\cite{Stefanini:2008} as follows:
\begin{definition}
Given two fuzzy numbers $u, v \in \fn$, the generalized Hukuhara difference (gH-difference for short) is
the fuzzy number $w$, if it exists, such that
\[
u\circleddash_{gH}v =w \Leftrightarrow\left\{ \begin{array}{c}
(i) \ u=v+w \\
\text{or}\\
(ii) \ v=u+(-1)w.
\end{array}\right.
\]
 It is easy to show that (i) and (ii) are both valid if and only if $w$ is a crisp number.
 If $u\circleddash_{gH}v$ exists, then its $\alpha$-cuts are given by
 \[
 [u\circleddash_{gH}v]_{\alpha} =[u]_{\alpha}\circleddash_{gH}[v]_{\alpha}
 \]
\end{definition}

Unfortunately, this substraction does not always define a fuzzy number.
The conditions to guarantee its existence are given in the next proposition

\begin{proposition}[\cite{Stefanini:2010}]\label{prop:existencegH}
Let $u,v\in \fn$ be two fuzzy numbers with $\alpha$-cuts given by $[u]_{\alpha}$ and
$[v]_{\alpha}$, respectively. The gH-difference $u\circleddash _{gH} v\in \fn$ if, and only if, one of the two conditions is satisfied
\begin{eqnarray}
(1) &\begin{cases}
len([u]_{\alpha})\geq len([v]_{\alpha}), \text{ for all $\alpha\in [0,1]$}, \\
u_\alpha^{-} - v_\alpha^{-} \text{ is increasing with respect to $\alpha$},\\
u_\alpha^{+} - v_\alpha^{+} \text{ is decreasing with respect to $\alpha$}.\\
\end{cases}\\
 \notag or& \\
(2) &\begin{cases}
len([u]_{\alpha})\leq len([v]_{\alpha}), \text{ for all $\alpha\in [0,1]$}, \\
u_\alpha^{+} - v_\alpha^{+} \text{ is increasing with respect to $\alpha$},\\
u_\alpha^{-} - v_\alpha^{-} \text{ is decreasing with respect to $\alpha$}.\\
\end{cases}\\ \notag
\end{eqnarray}
\end{proposition}

Let $f\in C([0,1],\fn)$. Take now the interval-valued function $f_{\alpha}:[0,1]\rightarrow \mathcal{K_C}$, defined as
$f_{\alpha}(x)=[f(x)]_{\alpha}$,  $\alpha\in[0,1]$.
So, for such functions $f\in C([0,1],\fn)$ we get a family of interval-valued functions $f_{\alpha}$ where our   Theorem~\ref{th:IV-degree1} in the previous section could be applied. Then we will need to retrieve the results back to our fuzzy-valued function $f$. Nevertheless,  as it was pointed out above,
a collection of closed intervals like
$$[u]_{\alpha}\circleddash_{gH}[v]_{\alpha}$$
can fail to define a proper fuzzy number.
For instance, as shown in \cite{BedeStef:2013}, given two triangular fuzzy numbers $u=\langle 12,15,19\rangle$ and $v=\langle 5,9,11\rangle$, $u\circleddash_{gH}v$ does not exist as $u^-_{\alpha}-v^-_{\alpha}=-6+\alpha$ and $u^+_{\alpha}-v^+_{\alpha}=-8+2\alpha$ are both increasing with respect to $\alpha\in[0,1]$, which fails to meet the conditions in Proposition~\ref{prop:existencegH}.
We will study that problem later. For now, we will add some restrictions to our continuous fuzzy-valued functions to prevent this from happening.

Firstly, we study how  the modulus of continuity of $f\in C([0,1],\fn)$ is related to the modulus of continuity of each level function $f_{\alpha}$.
We define the modulus of continuity of $f$ as
\[
\omega^{\mathcal{F}}(f,\delta):=\sup\{ d_{\infty}(f(x),f(y))\ :\ x,y\in[0,1],\  |x-y|<\delta\}.
\]
%and we recall that  the  modulus of continuity for each $f_{\alpha}$ is defined by
%\[
%\omega^{\mathcal{K}}(f_{\alpha},\delta)=\sup\{ d_H(f_{\alpha}(x),f_{\alpha}(y))\ :\ x,y\in[0,1],\  |x-y|<\delta\}.
%\]

\begin{lemma}\label{levels}
Let $f\in C([0,1],\fn)$. Then, for each  $\alpha\in[0,1]$, we have
\[
\omega^{\mathcal{K}}(f_{\alpha},\delta)\leq \omega^{\mathcal{F}}(f,\delta)
\]
\end{lemma}
\begin{proof}
Given any $\alpha\in[0,1]$, we have
\[
\begin{split}
\omega^{\mathcal{K}}(f_{\alpha},\delta)&=\sup\{ d_H(f_{\alpha}(x),f_{\alpha}(y))\ :\ x,y\in[0,1],\ |x-y|<\delta\}\\[2ex]
&\leq \sup\Big\{ \sup_{\alpha \in [0,1]} \{d_H(f_{\alpha}(x),f_{\alpha}(y))\}\ :\ x,y\in[0,1],\  |x-y|<\delta\Big\}\\
&=\sup\{ d_{\infty}(f(x),f(y))\ : \ x,y\in[0,1],\ |x-y|<\delta\}\\[2ex]
&=\omega^{\mathcal{F}}(f,\delta).
\end{split}
\]
\end{proof}

%In \cite{FSS:17} it is shown that the family
%
%$$\mathfrak{U}= \left\{\sum_{i=1}^{n} \psi_i(x)u_i\ :\ \psi_i \in C(K,[0,1]); \ u_i\in\fn,\ n=1,2,\ldots \right\}$$
%is dense in $C(K,\fn)$ with the supremum metric $D_{\infty}$, where $K$ is a compact metric space.
%Moreover,  it was proved that a subfamily is dense, since they got  the coefficient functions $\psi_i$ satisfying $\sum_{i=1}^{n} \psi_i =1$.

%When we consider the particular case of $K=[a,b]$, we can give a Jackson-type theorem but in order to get that we lose the condition $\sum_{i=1}^{n} \psi =1$, although we still can maintain
%$\sum_{i=1}^{n} \psi \leq 1$.

In order to apply Theorem~\ref{th:IV-degree1} we need to work with functions $f\in C([0,1],\fn)$
with the property that  $f(x)\circleddash _{gH}f(y)$ exists. Some examples of such functions are provided after Theorem~\ref{th:gH-fuzzy degree}

Let  us consider the family
$$\mathfrak{F}_n=\left\{\sum_{i=0}^{n} \psi_i(x)u_i\ :\ \psi_i \in C([0,1],[0,1]); \ u_i\in\fn,\ i=0,1,2,\ldots n\right\}.$$

Define the  approximation error between a continuous function $f\in C([0,1],\fn)$ and $\mathfrak{F}_n$ by
$$
E_{n,f}:=\inf_{g\in \mathfrak{F}_n}D(f,g).
$$

Let us recall that, given $u,v\in\fn$, the fuzzy inclusion $u\subseteq v$, means $[u]_{\alpha}\subseteq[v]_{\alpha}$, for all $\alpha\in [0,1]$. It is worth remarking that $u\subseteq v$ implies
$len([u]_{\alpha})\leq len([v]_{\alpha})$, for all $\alpha$.

The following theorem provides a Jackson-type approximation for such class of continuous functions.

\begin{theorem}\label{th:gH-fuzzy degree}
Let $f\in C([0,1],\fn)$ satisfying the following conditions:
\begin{enumerate}
\item There always exists $f(x)\circleddash _{gH}f(y)$, for all $0\leq x <y\leq 1$.
\item $f(y)\subseteq f(x)$, for all $x\leq y$.
\end{enumerate}
Then, for all $n\in\mathbb{N}$,
$$E_{n,f}\leq 2\omega^{\mathcal{F}}\left(f,\frac{1}{n}\right).$$
\end{theorem}
\begin{proof}
Fix $\varepsilon>0$ and $n\in\mathbb{N}$. Set $a_j=j/n$ and $u_j=f(a_j)$, $j=0,1,2,\ldots, n$.
%Being $f$ a continuous function it is bounded, so we can determine a constant $M>1$ such that
%$d_{\infty} (f(x),f(y))\leq 2M$, for all $x,y\in[0,1]$.
Take now $0<\delta<\frac{1}{2n} $.
% and
%$\epsilon'<\frac{\varepsilon}{2nM}$.
Condition (2) for the function $f$ implies that every function $f_{\alpha}$ satisfies equation~\eqref{eq:length1}. Then, we can
apply Proposition~\ref{prop:IV-degree2} to find $n+1$ continuous functions $\psi_j: [0,1]\rightarrow[0,1]$, such that, for each $\alpha\in[0,1]$, the interval-valued function

$$g_\alpha=f_{\alpha}(a_n)+\sum_{j=0}^{n-1} \psi_j(x) \Big(f_{\alpha}(a_j)\circleddash  f_{\alpha}(a_{j+1})\Big)$$
satisfies
\[
d_H(f_{\alpha}(x),g_{\alpha}(x))\leq 2\omega^{\mathcal{K}}(f_\alpha,1/n)+\varepsilon, \ \text{  for every $x\in [0,1]$}
\]
and, by Lemma~\ref{levels}, we get, for every $x\in [0,1]$,
\[
d_H(f_{\alpha}(x),g_{\alpha}(x))\leq 2\omega^{\mathcal{F}}(f,1/n)+\varepsilon.
\]
On the other hand, the condition (1)  on the function $f$ yields that the gH-difference $u_j\circleddash_{gH}u_{j+1}$ is well defined for all $j$. So we can define the fuzzy function $g\in\mathfrak{F}_n$ by
$$g(x):=u_n+\sum_{j=0}^{n-1} \psi_j(x) (u_j\circleddash _{gH} u_{j+1})$$
and then, for every $x\in[0,1]$, since $[u_j]_{\alpha}=f_{\alpha}(a_j)$, we have  $[g(x)]_{\alpha}=g_{\alpha}(x)$ and
\begin{align*}
d_{\infty}(f(x),g(x))&=\sup_{\alpha\in[0,1]} d_H([f(x)]_{\alpha},[g(x)]_{\alpha})\\[1ex]
&=\sup_{\alpha\in[0,1]} d_H(f_{\alpha}(x),g_{\alpha}(x))\leq 2\omega^{\mathcal{F}}(f,1/n)+\varepsilon,
\end{align*}
which means $D(f,g)\leq 2\omega^{\mathcal{F}}(f,1/n)+\varepsilon$. Finally
$$
E_{n,f}:=\inf_{h\in \mathfrak{F}_n}D(f,h)\leq D_{\infty}(f,g)\leq 2\omega^{\mathcal{F}}(f,1/n)+\varepsilon
$$
and, being true for all $\varepsilon>0$, we get the conclusion.

\end{proof}

\begin{example}\normalfont
Let $g$ be a positive real-valued decreasing continuous function defined on $[0,1]$.
Let $f\in C([0,1],\fn)$ be defined as $f(x):=g(x)\cdot u$ being $u$ a fixed fuzzy number. Then, for every $0\leq x<y\leq 1$, there always exists $f(x)\circleddash _{gH} f(y)$.
Since $g(x)$ is decreasing, $len([f(x)]_{\alpha})\geq len([f(y)]_{\alpha})$, for every $\alpha\in [0,1]$. On the other hand,
if $[u]_{\alpha}=[u_{\alpha}^{-}, u_{\alpha}^{+}]$, then
\begin{eqnarray*}
w_{\alpha}^{-} &=& g(x)u_{\alpha}^{-}- g(y)u_{\alpha}^{-}= (g(x)-g(y))u_{\alpha}^{-}, \\[2ex]
 w_{\alpha}^{+}&=& g(x)u_{\alpha}^{+}- g(y)u_{\alpha}^{+}= (g(x)-g(y))u_{\alpha}^{+}
\end{eqnarray*}
and $w_{\alpha}^{-}$ is increasing and $w_{\alpha}^{+}$ is decreasing with respect to $\alpha$, which means that   the condition $(1)$ in Proposition~\ref{prop:existencegH}  is fulfilled and  then,
it always defines a fuzzy number $f(x)\circleddash _{gH}f(y)\in\fn$ with $\alpha$-cuts
\[
[f(x)\circleddash _{gH} f(y)]_{\alpha}=[w_{\alpha}^{-},w_{\alpha}^{+}]
=[f(x)]_{\alpha}\circleddash _{gH}[f(y)]_{\alpha}
\]
It is worth mentioning  that if the function $g$ is increasing, then the function $f(x):=g(x)\cdot u$ also satisfies that
 $f(x)\circleddash _{gH}f(y)\in\fn$, fulfilling now condition (2) in Proposition~\ref{prop:existencegH}.

 \end{example}

%A basic result in fuzzy number theory is the so-called
%representation theorem of fuzzy numbers on $\mathbb{R}$, which was
%provided by Goetschel and Voxman in \cite{GV}:

Using Proposition~\ref{prop:IV-degree3} we can show, similarly,
\begin{theorem}\label{th:gH-fuzzy degree2}
Let $f\in C([0,1],\fn)$ satisfying the following conditions:
\begin{enumerate}
\item There always exists $f(y)\circleddash _{gH}f(x)$, for all $0\leq x <y\leq 1$.
\item $f(y)\supseteq f(x)$, for all $x\leq y$.
\end{enumerate}
Then, for all $n\in\mathbb{N}$,
$$E_{n,f}\leq 2\omega^{\mathcal{F}}\left(f,\frac{1}{n}\right).$$
\end{theorem}

As mentioned above, the $gH$-difference of fuzzy numbers does not always exist. To overcome this problem, Stefanini proposed a new difference between fuzzy numbers that always exists:
\begin{definition}(\cite{Stefanini:2010},\cite{BedeStef:2013},\cite{GomesBarros:15})
Given two fuzzy numbers $u$, $v$, the generalized difference (g-difference, for short) of $u$ and $v$ is given by their level sets as
\[
[u\circleddash_g v]_{\alpha}=cl\Big(conv\bigcup_{\beta\geq\alpha} [u]_{\beta}\circleddash[v]_{\beta}\Big),
\]
where the gH-difference $\circleddash$ is used with interval operands $[u]_{\beta}$ and $[v]_{\beta}$.
\end{definition}
It turns out that the $\alpha$-cuts can be calculated as
\[
[u\circleddash_g v]_{\alpha} =
[\inf_{\beta\geq\alpha} \min\{ u_{\beta}^{-}-v_{\beta}^{-}, u_{\beta}^{+}-v_{\beta}^{+}\},
\sup_{\beta\geq\alpha} \max\{ u_{\beta}^{-}-v_{\beta}^{-}, u_{\beta}^{+}-v_{\beta}^{+}\}]
\]

Although, the gH-difference of intervals $[u]_{\alpha}\circleddash [v]_{\alpha}$ could not coincide with
$[u\circleddash_g v]_{\alpha}$, except when $u\circleddash_{gH} v$ exists, there is still some connection for the Hausdorff metric.

\begin{lemma}\label{lem:g-gH-cuts}
Let $u,v\in \fn$ be two fuzzy numbers with their $\alpha$-cuts satisfying that $[v]_{\alpha}\subseteq [u]_{\alpha}$, for all $\alpha\in [0,1]$. Then,
\[
d_H \Big([u]_{\alpha}\circleddash  [v]_{\alpha},\ [u\circleddash _g v]_{\alpha}\Big)\leq
\sup_{\beta\geq \alpha} d_H \Big([u]_{\alpha}\circleddash  [v]_{\alpha},\ [u]_{\beta}\circleddash  [v]_{\beta}\Big).
\]

\end{lemma}
\begin{proof}
The condition $[v]_{\alpha}\subseteq [u]_{\alpha}$
gives $u_{\beta}^{-}-v_{\beta}^{-}\leq u_{\beta}^{+}-v_{\beta}^{+}$, for all $\beta$, which
implies that
\begin{eqnarray*}
[u\circleddash _g v]_{\alpha} &=&[\inf_{\beta\geq \alpha}\{u_{\beta}^{-} - v_{\beta}^{-}\},
\sup_{\beta\geq \alpha}\{u_{\beta}^{+} - v_{\beta}^{+}\}];\\[2ex]
[u]_{\alpha}\circleddash  [v]_{\alpha}&=& [u_{\alpha}^{-} - v_{\alpha}^{-},
\ u_{\alpha}^{+} - v_{\alpha}^{+}]
\end{eqnarray*}
so
\[
\begin{split}
d_H \Big([u]_{\alpha}\circleddash  [v]_{\alpha},\ [u\circleddash _g v]_{\alpha}\Big) &=
\max\{\Big| u_{\alpha}^{-} - v_{\alpha}^{-}- \inf_{\beta\geq \alpha}\{u_{\beta}^{-} - v_{\beta}^{-}\}\Big|,
\ \Big| u_{\alpha}^{+} - v_{\alpha}^{+}-\sup_{\beta\geq \alpha}\{u_{\beta}^{+} - v_{\beta}^{+}\}\Big| \}.
\end{split}
\]
Assume that
\[
d_H \Big([u]_{\alpha}\circleddash  [v]_{\alpha},\ [u\circleddash _g v]_{\alpha}\Big)=
\Big| u_{\alpha}^{-} - v_{\alpha}^{-}- \inf_{\beta\geq \alpha}\{u_{\beta}^{-} - v_{\beta}^{-}\}\Big|.
\]
Then, removing the absolute value, we have
\[
\begin{split}
d_H \Big([u]_{\alpha}\circleddash  [v]_{\alpha},\ [u\circleddash _g v]_{\alpha}\Big)&=
u_{\alpha}^{-} - v_{\alpha}^{-}- \inf_{\beta\geq \alpha}\{u_{\beta}^{-} - v_{\beta}^{-}\}\\[2ex]
&=\sup_{\beta\geq \alpha}\{u_{\alpha}^{-} - v_{\alpha}^{-} - (u_{\beta}^{-} - v_{\beta}^{-})\} \\[2ex]
&\leq \sup_{\beta\geq \alpha} d_H \Big([u]_{\alpha}\circleddash  [v]_{\alpha},\ [u]_{\beta}\circleddash  [v]_{\beta}\Big).
\end{split}
\]
On the other case, if
\[
d_H \Big([u]_{\alpha}\circleddash  [v]_{\alpha},\ [u\circleddash _g v]_{\alpha}\Big)=
\Big| u_{\alpha}^{+} - v_{\alpha}^{+}- \sup_{\beta\geq \alpha}\{u_{\beta}^{+} - v_{\beta}^{+}\}\Big|
\]
then,
\[
\begin{split}
d_H \Big([u]_{\alpha}\circleddash  [v]_{\alpha},\ [u\circleddash _g v]_{\alpha}\Big)&=
\sup_{\beta\geq \alpha}\{u_{\beta}^{+} - v_{\beta}^{+}\} -(u_{\alpha}^{+} - v_{\alpha}^{+})\\[2ex]
&=\sup_{\beta\geq \alpha}\{u_{\beta}^{+} - v_{\beta}^{+}-(u_{\alpha}^{+} - v_{\alpha}^{+})\}\\[2ex]
&\leq \sup_{\beta\geq \alpha} d_H \Big([u]_{\alpha}\circleddash  [v]_{\alpha},\ [u]_{\beta}\circleddash  [v]_{\beta}\Big)
\end{split}
\]
and we are done.

\end{proof}

Unfortunately, the generalized difference of two fuzzy numbers has no cancelation law and the technique used in Theorem~\ref{th:gH-fuzzy degree}
is not applicable.
Using  Lemma~\ref{lem:g-gH-cuts} we are still able to provide a
Jackson-type approximation for continuous fuzzy-valued functions, avoiding one of the restrictions in Theorem~\ref{th:gH-fuzzy degree}, although with a nonconstant bound.

\begin{theorem}\label{th:main-fuzzy-degree}
Let $f\in C([0,1],\fn)$ satisfying
$$f(y)\subseteq f(x), \text{ for all $x\leq y$.}$$

Then, for all $n\in\mathbb{N}$,
$$E_{n,f}\leq (2n+2)\omega^{\mathcal{F}}\left(f,\frac{1}{n}\right).$$
\end{theorem}
\begin{proof}
Proceed as in Theorem~\ref{th:gH-fuzzy degree}. Fix $\varepsilon>0$ and $n\in\mathbb{N}$. Set $a_j=j/n$ and $u_j=f(a_j)$, $j=0,1,2,\ldots, n$ and find $n+1$ continuous functions $\psi_j: [0,1]\rightarrow[0,1]$ satisfying conditions in Proposition~\ref{prop:IV-degree2} for
$0<\delta<\frac{1}{2n} $.
% and
%$\epsilon'<\frac{\varepsilon}{2nM}$, being $M$ a constant with $d_{\infty} (f(x),f(y))\leq 2M$, for all $x,y\in[0,1]$.

\medskip
We already know
 that, for each $\alpha\in[0,1]$, the interval-valued function

$$g_\alpha=f_{\alpha}(a_n)+\sum_{j=0}^{n-1} \psi_j(x) \Big(f_{\alpha}(a_j)\circleddash  f_{\alpha}(a_{j+1})\Big)$$
satisfies, for every $x\in [0,1]$,
\[
d_H(f_{\alpha}(x),g_{\alpha}(x))\leq 2\omega^{\mathcal{F}}(f,1/n)+\varepsilon
\]
where $f_{\alpha}(x)=[f(x)]_{\alpha}$.

%\medskip
%{\color{red} MAL
%
%Let us consider $A=\max_{x\in[0,1]} \sum_{j=0}^{n-1}\psi_j(x)>0$. If $A\leq 1$, then take $B=1$.  Otherwise, take $B=\min\{\frac{1}{A},1-\frac{1}{A}\}$. It is worth to remark that both $B$ and $1-B$ are less than $\frac{1}{A}$.
%It is apparent that, for every $x\in[0,1]$,
%$$\frac{1}{A}\sum_{j=0}^{n-1} \psi_j(x) \leq 1$$
%so
%$$\sum_{j=0}^{n-1} B\psi_j(x)\leq 1 \text{ and }\sum_{j=0}^{n-1} (1-B)\psi_j(x)\leq 1 $$}
%

\medskip
Let us consider the continuous fuzzy-valued function in the family $\mathcal{F}_n$ defined by
$$h(x):=u_n+\sum_{j=0}^{n-1} \psi_j(x) (u_j\circleddash_g u_{j+1}).$$

Fix $x\in[0,1]$ and $\alpha\in[0,1]$. Since
\[
\begin{split}
d_H\left([f(x)]_{\alpha}, [h(x)]_{\alpha}\right)&\leq d_H([f(x)]_{\alpha}, g_{\alpha}(x))
+d_H(g_{\alpha}(x), [h(x)]_{\alpha})\\
&\leq 2\omega^{\mathcal{F}}(f,1/n)+\varepsilon +d_H(g_{\alpha}(x), [h(x)]_{\alpha}),\\
\end{split}
\]
we need to compute
$d_H(g_{\alpha}(x), [h(x)]_{\alpha})$. Since, by Lemma~\ref{lem:g-gH-cuts}
\[
\begin{split}
d_H\Big([u_j]_{\alpha}\circleddash [u_{j+1}]_{\alpha}),
[u_j\circleddash_g u_{j+1}]_{\alpha}\Big)&\leq
\sup_{\beta\geq \alpha} d_H \Big([u_j]_{\alpha}\circleddash  [u_{j+1}]_{\alpha},\ [u_j]_{\beta}\circleddash  [u_{j+1}]_{\beta}\Big)\\
&\leq \sup_{\beta\geq \alpha} \Big[
d_H \Big([u_j]_{\alpha}\circleddash  [u_{j+1}]_{\alpha}, \mathbf{0}\Big) +
d_H\Big( [u_j]_{\beta}\circleddash  [u_{j+1}]_{\beta},\mathbf{0}\Big)\Big]\\
&\leq \omega^{\mathcal{F}}(f,1/n) +\omega^{\mathcal{F}}(f,1/n),
\end{split}
\]
we get
\[
\begin{split}
d_H\left(g_{\alpha}(x), [h(x)]_{\alpha}\right)&\leq
\sum_{j=0}^{n-1}d_H\Big(\psi_j(x)([u_j]_{\alpha}\circleddash [u_{j+1}]_{\alpha}),
\psi_j(x) [u_j\circleddash_g u_{j+1}]_{\alpha}\Big)\\
&\leq
\sum_{j=0}^{n-1}\psi_j(x)d_H\Big([u_j]_{\alpha}\circleddash [u_{j+1}]_{\alpha},
[u_j\circleddash_g u_{j+1}]_{\alpha}\Big)\\
&\leq 2\omega^{\mathcal{F}}(f,1/n)\sum_{j=0}^{n-1}\psi_j(x)\\
&\leq 2n\omega^{\mathcal{F}}(f,1/n).
\end{split}
\]
So, we have just shown that
\[
d_H\left([f(x)]_{\alpha}, [h(x)]_{\alpha}\right)
\leq 2\omega^{\mathcal{F}}(f,1/n)+\varepsilon +2n\omega^{\mathcal{F}}(f,1/n)
\]
for  every $x$ and $\alpha$, which means
\[
D(f,g)\leq (2n+2)\omega^{\mathcal{F}}(f,1/n)+\varepsilon.
\]
Then,
\[
\inf_{h\in \mathfrak{F}_n}D(f,h)\leq D(f,g)\leq (2n+2)\omega^{\mathcal{F}}(f,1/n)+\varepsilon
\]
and, being true for all $\varepsilon>0$, we finally get the conclusion
\[
E_{n,f}\leq (2n+2)\omega^{\mathcal{F}}(f,1/n).
\]
\end{proof}
The problem with the nonconstant bound seems to be that the sum of the coefficients $\sum \psi_k$ is not bounded. So, in order
to improve this result, we need to find a family of functions $\phi_k$ well behaved for the modulus of continuity  and with the property that $\sum_{k=0}^{n} \phi_k=1$. Thus, using the technique of Prolla (\cite{prola}), we are able to provide a Jackson type approximation, with a slight worse bound than the one in Theorem~\ref{th:gH-fuzzy degree}, but with the advantage that
it can be applied to any continuous fuzzy-valued function defined on $[0,1]$.

\begin{proposition}\label{th:sum2}
There exist $n+1$ continuous functions $\phi_j:[0,1]\rightarrow[0,1]$ such that,  for every $x\in[0,1]$,
$$\sum_{j=0}^{n}\phi_j(x)=1,$$
and with the property that, at each $x\in[0,1]$,  there exist, at most, two nonzero terms  in the sum.
\end{proposition}
\begin{proof}
Take $a_j=\frac{j}{n}$, $j=0,1,\ldots, n$ and let $0<\delta<\frac{1}{2n}$. Define the continuous functions $f_j:[0,1]\rightarrow[0,1]$, $j=0,1,\ldots, n$, as follows

$$
\begin{aligned}
f_0(x)&:=\begin{cases}
\frac{1}{\delta}(\delta-x), & \text{if $0 \leq x< \delta$}\\
0, & \text{if $\delta\leq x\leq 1$}
\end{cases}\\
f_j(x)&:=\begin{cases}
1, & \text{if $0\leq x\leq a_j-\delta$}\\
\frac{1}{2\delta}(a_j+\delta-x), & \text{if $a_j-\delta < x< a_j+\delta$}\\
0, & \text{if $a_j+\delta\leq x\leq 1$}
\end{cases}\\
f_n(x)&:=\begin{cases}
1, & \text{if $0\leq x\leq a_n-\delta$}\\
\frac{1}{\delta}(a_n-x), & \text{if $a_n-\delta < x\leq 1$}\\
\end{cases}
\end{aligned}
$$
that is, $f_j$ takes value $1$ before $a_j-\delta$, takes value $0$ after $a_j+\delta$ and it is linear in the gap between.
Now, let us define,
\[
\begin{aligned}
\phi_0&=f_0,\\
\phi_1&=(1-f_0)f_1,\\
\phi_2&=(1-f_0)(1-f_1)f_2,\\
&\vdots\\
\phi_{n-1}&=(1-f_0)(1-f_1)\cdots (1-f_{n-2})f_{n-1}.
\end{aligned}
\]
We claim that
$$\phi_0+\phi_1+\phi_2+\ldots+\phi_{n-1}=
1-(1-f_0)(1-f_1)\cdots(1-f_{n-1}).$$

Then, taking $\phi_n=(1-f_0)(1-f_1)\cdots(1-f_{n-1})$ we have proved that
$$\sum_{j=0}^{n} \phi_j(x)=1,\quad\text{for all $x\in [0,1]$}.$$

We have partitioned $[0,1]$ into two disjoint sets, say,
$$
\begin{aligned}
V&:=[0,\delta)\cup\Big(\bigcup_{j=1}^{n-1}(a_j-\delta,a_j+\delta)\Big)\cup(a_n-\delta,1];\\
W&:=[0,1]\setminus V.
\end{aligned}
$$
We will study what happens to each $\phi_j(x)$ depending on where $x$ is located.

Assume that $x\in V$. Then $x\in(a_k-\delta,a_k+\delta)$ for some $k\in\{1,\ldots,n-1\}$ or $x\in [0,\delta)$ or $x\in(a_n-\delta,1]$.

If $x\in [0,\delta)$, then $f_1(x)=1$ and $\phi_j(x)=0$, for all $j>1$.
We  have only two nonzero summands, $\phi_0(x)$ and $\phi_1(x)$.
If $x\in(a_n-\delta,1]$ then $f_j(x)=0$ and therefore $\phi_j(x)=0$, for all $0\leq j\leq n-1$. So we have only one nonzero term $\phi_n(x)=1$.

If $x\in(a_k-\delta,a_k+\delta)$, $1\leq k\leq n-1$, then
$$x> a_{k-1}-\delta>a_{k-2}-\delta>\ldots> \delta,$$
implying that $f_j(x)=0$ and $\phi_j(x)=0$ as well,  for $0\leq j\leq k-1$.
On the other hand,
$$x<a_{k+1}-\delta<a_{k+2}-\delta<\ldots<a_n-\delta,$$
which yields $f_j(x)=1$ for $j\geq k+1$ and  then $\phi_j(x)=0$, for $j>k+1$.
Again, that gives us only two nonzero terms, $\phi_k(x)$ and $\phi_{k+1}(x)$.

We study now the case $x\in W$. This implies that $x\in[a_{k-1}+\delta, a_k-\delta]$, for some $1\leq k\leq n$.
Then
$$x\leq a_{k}-\delta<a_{k+1}-\delta<\ldots<a_n-\delta,$$
which means $f_j(x)=1$ for $j\geq k$ and  then $\phi_j(x)=0$, for $j\geq k+1$.
On the other hand,
$$x\geq a_{k-1}+\delta>a_{k-2}+\delta>\ldots> \delta,$$
implying that $f_j(x)=0$ and also $\phi_j(x)=0$,  for $0\leq j\leq k-1$.
Again we have only  one nonzero term, $\phi_k(x)=1$.
\end{proof}

The modulus of continuity has the following properties (\cite{Gal:2000}).
\begin{proposition}\label{prop:modulus1}
Let $f$ be in $C([0,1],\fn)$. Then,
\begin{enumerate}
\item $\omega^{\mathcal{F}}(f,\delta_1+\delta_2)\leq \omega^{\mathcal{F}}(f,\delta_1)+\omega^{\mathcal{F}}(f,\delta_2)$, $\delta_1,\delta_2>0$.
\item $\omega^{\mathcal{F}}(f,n\delta)\leq n\;\omega^{\mathcal{F}}(f,\delta)$, $\delta>0$, $n\in\mathbb{N}$.
\item $\omega^{\mathcal{F}}(f,\lambda\delta)\leq (\lambda+1)\;\omega^{\mathcal{F}}(f,\delta)$, $\delta>0$, $\lambda>0$.
\end{enumerate}

\end{proposition}

The functions $\phi_k(x)$ in Proposition~\ref{th:sum2} are piecewise linear and, in fact, they have a trapezoidal graphic. For each $0<\delta<\frac{1}{2n}$, we define:
$$T_{n,\delta}^{f}(x):=\sum_{k=0}^{n} \phi_k(x)\cdot f\left(\frac{j}{n}\right).$$
Let $\mathfrak{T}_{n,\delta}$ be the set of all such functions.
Now we can prove the following theorem.
\begin{theorem}\label{th:degree-approx-fuzzy}
Let $f\in C([0,1],\fn)$.
Then, for all $n\in\mathbb{N}$ and $0<\delta<\frac{1}{2n}$,
$$D\left(f,T_{n,\delta}^{f}\right)\leq 3\omega^{\mathcal{F}}\left(f,\frac{1}{n}\right).$$
\end{theorem}
\begin{proof}
Set $a_j=\frac{j}{n}$, $j=0,1,\ldots,n$. Set $u_j=f(a_j)$ and take $0<\delta<\frac{1}{2n}$. Let us take  $g(x):=T_{n,\delta}^{f}(x)$.

Then, for each $x\in[0,1]$ and $\alpha\in[0,1]$,
\[
\begin{split}
d_{H}([f(x)]_{\alpha}, [g(x)]_{\alpha}) &=
d_{H}\Big(\sum_{j=0}^{n}\phi_j(x)[f(x)]_{\alpha}, \sum_{j=0}^{n}\phi_j(x)[u_j]_{\alpha}\Big) \\
&\leq \sum_{j=0}^{n}\phi_j(x)d_{H}\Big([f(x)]_{\alpha},[u_j]_{\alpha}\Big).
\end{split}
\]
Let us consider
$$
\begin{aligned}
V&:=[0,\delta)\cup\Big(\bigcup_{j=1}^{n-1}(a_j-\delta,a_j+\delta)\Big)\cup(a_n-\delta,1],\\
W&:=[0,1]\setminus V.
\end{aligned}
$$
Now, if $x\in V$, then we have two nonzero terms, $\phi_k(x)$ and $\phi_{k+1}(x)$ for $k=0,1,\ldots, n-1$ when $x\in(a_k-\delta,a_k+\delta)$ or $x\in[0,\delta)$ and we have only one $\phi_n(x)=1$, when $x\in(a_n-\delta,1]$. Then, in the worst case and bearing in mind that the position of $x$ means $|x-a_k|<\frac{1}{n}$ and $|x-a_{k+1}|<\frac{2}{n}$, we have
$$
\begin{aligned}
d_{H}\Big([f(x)]_{\alpha},[u_k]_{\alpha}\Big)&\leq \omega^{\mathcal{F}}\left(f,\frac{1}{n}\right)\quad\text{and}\\
d_{H}\Big([f(x)]_{\alpha},[u_{k+1}]_{\alpha}\Big)&\leq \omega^{\mathcal{F}}\left(f,\frac{2}{n}\right).\\
\end{aligned}
$$
Then, we obtain
\[
\begin{split}
d_{H}([f(x)]_{\alpha}, [g(x)]_{\alpha}) &\leq
 \phi_k(x)d_{H}\Big([f(x)]_{\alpha},[u_k]_{\alpha}\Big)+
 \phi_{k+1}(x)d_{H}\Big([f(x)]_{\alpha},[u_{k+1}]_{\alpha}\Big)\\
 &\leq \omega^{\mathcal{F}}\left(f,\frac{1}{n}\right) + \omega^{\mathcal{F}}\left(f,\frac{2}{n}\right) \leq \omega^{\mathcal{F}}\left(f,\frac{1}{n}\right)+2\omega^{\mathcal{F}}\left(f,\frac{1}{n}\right)\\
 &=3\omega^{\mathcal{F}}\left(f,\frac{1}{n}\right)
\end{split}
\]
Now, when $x\in W$, there exists only one nonzero term $\phi_k(x)$. So,
\[
\begin{split}
d_{H}([f(x)]_{\alpha}, [g(x)]_{\alpha}) &\leq
 \phi_k(x)d_{H}\Big([f(x)]_{\alpha},[u_k]_{\alpha}\Big)\\
 &\leq \omega^{\mathcal{F}}\left(f,\frac{1}{n}\right)
\leq 3\omega^{\mathcal{F}}\left(f,\frac{1}{n}\right).
\end{split}
\]
So, we have just shown that,
\[
d_{H}([f(x)]_{\alpha}, [g(x)]_{\alpha}) \leq  3\omega^{\mathcal{F}}\left(f,\frac{1}{n}\right),
\]
for every $x\in[0,1]$ and for every $\alpha$-cut, which means that
$$D(f,T_{n,\delta}^{f})\leq 3\omega^{\mathcal{F}}\left(f,\frac{1}{n}\right).$$
%and then
%$$E_{n,f}\leq D(f,g)\leq 3\omega^{\mathcal{F}}\left(f,\frac{1}{n}\right).$$
\end{proof}

\end{document}